\documentclass[12pt]{article}

\usepackage{amsfonts}

\usepackage{amscd}
\usepackage{amsthm}
\usepackage{indentfirst}

\usepackage[all]{xy}

\usepackage{amssymb, amsmath, amsthm, amsgen, amstext, amsbsy, amsopn}
\usepackage{epsfig}
\usepackage{epic,eepic}

\newtheorem{thm}{Theorem}[section]
\newtheorem{lemma}[thm]{Lemma}
\newtheorem{prop}[thm]{Proposition}
\newtheorem{cor}[thm]{Corollary}

\newtheorem{rem}[thm]{Remark}

\newcommand{\R}{{\mathbb{R}}}

\newcommand{\Z}{{\mathbb{Z}}}
\newcommand{\N}{{\mathbb{N}}}
\newcommand{\C}{{\mathbb{C}}}

\newcommand{\cH}{{\mathcal{H}}}

\newcommand{\cJ}{{\mathcal{J}}}
\newcommand{\cL}{{\mathcal{L}}}

\newcommand{\cP}{{\mathcal{P}}}
\newcommand{\cQ}{{\mathcal{Q}}}
\newcommand{\cR}{{\mathcal{R}}}
\newcommand{\cS}{{\mathcal{S}}}

\newcommand{\cU}{{\mathcal{U}}}
\newcommand{\cV}{{\mathcal{V}}}
\newcommand{\cW}{{\mathcal{W}}}

\newcommand{\hcQ}{{\widehat{\mathcal{Q}}}}
\def\id{{1\hskip-2.5pt{\rm l}}}

\newcommand{\gog}{{\mathfrak{g}}}
\newcommand{\goh}{{\mathfrak{h}}}
\newcommand{\gou}{{\mathfrak{u}}}
\newcommand{\spn}{{\mathfrak{sp}\, (2n,\R)}}
\newcommand{\gogl}{{\mathfrak{gl}\, (n,\R)}}

\newcommand{\gt}{\mathfrak}

\newcommand{\cC}{{\mathcal{C}}}

\newcommand{\tSp}{{\widetilde{Sp}}}

\newcommand{\diam}{{\hbox{\rm diam}\,}}

\newcommand{\ad}{{\hbox{\it ad}}}

\newcommand{\supp}{{\it supp\,}}

\newcommand{\Qed}{\hfill \qedsymbol \medskip}

\begin{document}

\title{Lie quasi-states \\}

\renewcommand{\thefootnote}{\alph{footnote}}

\author{\textsc Michael Entov$^{a}$,\ Leonid
Polterovich$^{b}$}

\footnotetext[1]{Partially supported by the Israel Science
Foundation grant $\#$ 881/06.} \footnotetext[2]{Partially
supported by the Israel Science Foundation grant $\#$ 509/07.}

\date{\today}

\maketitle

\begin{abstract}
\noindent Lie quasi-states on a real Lie algebra are functionals
which are linear on any abelian subalgebra. We show that on the
symplectic Lie algebra of rank at least 3 there is only one
continuous non-linear Lie quasi-state (up to a scalar factor,
modulo linear functionals). It is related to the asymptotic Maslov
index of paths of symplectic matrices.

\end{abstract}

\tableofcontents


\renewcommand{\thefootnote}{\arabic{footnote}}
\vfil \eject

\section{Introduction and main results}
\label{sec-intro}

\subsection{Lie quasi-states}

Let $W\subset \gog$ be a vector subspace of a finite-dimensional
Lie algebra $\gog$ over $\R$. A function $\zeta: W \to \R$ will be
called {\it quasi-linear} if:
$$
[x_1,x_2]=0\ \Longrightarrow\ \zeta (c_1 x_1 + c_2 x_2) = c_1 \zeta (x_1) + c_2 \zeta
(x_2)\
\forall c_1, c_2\in\R.
$$
A quasi-linear function on the whole Lie algebra $\gog$ will be
called a {\it Lie quasi-state}.

Continuous Lie quasi-states on $\gog$ form a vector space $\hcQ (\gog)$. Set
$\cQ (\gog) := \hcQ (\gog)/\gog^*$, where
$\gog^*$ is the dual space to $\gog$.  It can be viewed as the space of {\it
non-linear} continuous Lie quasi-states on $\gog$.

In the present paper we focus on Lie quasi-states on the
symplectic Lie algebra $\spn$, that is on the Lie algebra of the
group $Sp (2n,\R)$ of $2n\times 2n$ symplectic matrices. Our main
finding is that the notion of a continuous Lie quasi-state is
rigid in the following sense.

\begin{thm}
\label{thm-main} Let $\gog = \spn$, $n\geq 3$. Then $\dim \cQ
(\gog) = 1$.
\end{thm}

\medskip
\noindent As we shall see below, the generator of $\cQ (\gog)$
looks as follows: its value on a matrix $B \in \spn$ equals,
roughly speaking, to the asymptotic Maslov index of the path
$e^{tB}$ as $t \to \infty$.

\medskip

Let us discuss the assumptions of the theorem.  In the case $n=1$,
 $\dim \cQ (\mathfrak{sp}\, (2, \R)) =
+\infty$. Indeed, any two commuting matrices in $\mathfrak{sp}\, (2,
\R)$ differ by a scalar factor, and hence any odd homogeneous
function on $\mathfrak{sp}\, (2, \R)$ is a Lie quasi-state. The case
$n=2$ is so far absolutely open.

\medskip \noindent The next result shows that the continuity
assumption in Theorem~\ref{thm-main} is essential.

\begin{thm}
\label{thm-discont} The space of (not necessarily continuous) Lie
quasi-states on  $\spn$ which are bounded on a neighborhood of zero
is infinite-dimen\-si\-onal for all $n \geq 1$.
\end{thm}

\medskip
\noindent At the same time any Lie quasi-state which is {\it
differentiable} at $0$ is automatically linear since it is
homogeneous of degree $1$.

\subsection{Origins of Lie quasi-states}

The interest to the notion of Lie quasi-states is
three-fold.

\medskip
\noindent {\sc Lie quasi-states and quasi-morphisms on Lie
groups:} Recall that a {\it homogeneous quasi-morphism} on a group
$G$ is a function $\mu : G\to\R$ such that

\begin{itemize}

\item{} There exists $C>0$ so that $|\mu (xy) - \mu (x) - \mu
(y)|\leq C$ for all $x,y\in G$.

\item{} $\mu (x^k) = k\mu (x)$ for all $k\in\Z$, $x\in G$.

\end{itemize}

It is known that restriction of any homogeneous quasi-morphism to an
abelian subgroup is a genuine morphism, and that homogeneous
quasi-mor\-phisms are conjugation invariant (see e.g. \cite{Bav} for
introduction to quasi-mor\-phisms). Therefore, given a homogeneous
quasi-morphism $\mu$ on a Lie group $G$, its pull-back to the Lie
algebra $\gog$ by the exponential map, which we will call {\it the
directional derivative of $\mu$},
$$\zeta: \gog \to \R, \; a \mapsto \mu(\exp a)\;,$$
is an $Ad_G$-invariant Lie quasi-state. Clearly, $\zeta$ is
continuous whenever $\mu$ is continuous.

In fact, if $G$ is a simply connected real Lie group with a
Hermitian simple Lie algebra $\gog$ (see
Section~\ref{sec-Ad-inv-Lie-qs} for the definition; in particular,
$G=\tSp (2n,\R)$, the universal cover of $Sp (2n,\R)$, is such a Lie
group and its Lie algebra is $\spn$), then the space of homogeneous
quasi-morphisms on $G$ is 1-dimensional \cite{Shtern}, cf.
\cite{BS-H}, and these quasi-morphisms are continuous \cite{Shtern}.
We shall show that in such a case the space of $Ad_G$-invariant Lie
quasi-states on $\gog$ is also one-dimensional -- see
Section~\ref{sec-Ad-inv-Lie-qs}.

\medskip
\noindent {\sc Lie quasi-states and Gleason's theorem:} Gleason's
theorem \cite{Gleason} is one of the most famous and important
results in the mathematical formalism of quantum mechanics (see e.g.
\cite{Peres, Dvurecenskij}). In the finite-dimensional setting the
proof of Gleason's theorem yields the following result about Lie
quasi-states.

\begin{thm}[Gleason]
\label{thm-Gleason} Let $V$ be a finite-dimensional vector space
over $\R$ (respectively, over $\C$), equipped with a real
(respectively, Hermitian) inner product. Denote by $\cS(V)$ the
subspace of the self-adjoint operators (viewed as the subspace of
the Lie algebra of all operators on $V$). Let $\zeta: \cS(V)\to\R$
be a quasi-linear function which is bounded on a neighborhood of
zero in $\cS (V)$. Assume also that $\dim V\geq 3$.

Then $\zeta$ is linear and has the form $\zeta (A) = tr (HA)$ for
some $H\in \cS(V)$.

\end{thm}

\begin{cor}
\label{cor-Gleason} Any Lie quasi-state $\zeta$ on the Lie algebra
$\gou(n)$, $n\geq 3$, which is bounded on a neighborhood of zero, is
linear and has the form $\zeta (A) = tr (HA)$ for some $H\in
\gou(n)$.

\end{cor}

Indeed, $\gou (n) = i\cS (\C^n)$, where $\cS (\C^n)$ is the space of
Hermitian $n\times n$-matrices.

The statement of Theorem~\ref{thm-Gleason} is slightly different
from the original formulation in Gleason's paper \cite{Gleason}:
instead of boundedness of $\zeta$ near zero Gleason assumes that
$\zeta$ is non-negative on the set of non-negative self-adjoint
operators. To obtain this non-negativity condition from the
boundedness near zero one just needs to add to $\zeta$ a linear
function $A\mapsto Tr (cA)$ for a sufficiently large positive
$c\in\R$. After such a modification, $\zeta$ becomes positive on the
set of all orthogonal projectors, and hence (by the spectral
theorem) non-negative on all non-negative self-adjoint operators.

Our proof of Theorem~\ref{thm-main} uses Gleason's theorem.  Let
us mention that the most difficult and non-trivial part of the
proof of Gleason's Theorem~\ref{thm-Gleason} is to show that the
boundedness of $\zeta$ near zero implies its continuity -- the
latter yields (by basic representation theory) that $\zeta$ is
linear. Since in Theorem~\ref{thm-main} we assume that $\zeta$ is
continuous, our proof of the theorem does not use the difficult
part of Gleason's proof (and there is no analogue of this part in
our proof due to Theorem~\ref{thm-discont}).

\medskip
\noindent{\sc Lie quasi-states in symplectic topology:}  As the
third point of interest in Lie quasi-states, we note that such
functionals on the infinite-dimensional Poisson-Lie algebra of
Hamiltonian functions on a symplectic manifold appeared recently in
symplectic topology and Hamiltonian dynamics before they were
properly studied in the finite-dimensional setting. We refer the
reader to \cite{EP-qst,EPZ,EPZ-physics} for various aspects of this
development.

\subsection{Maslov quasi-state on $\spn$}

The Lie quasi-state $\zeta_M$ generating the 1-dimensional space
$\cQ (\spn)$, $n\geq 3$, comes from the Maslov index of paths of symplectic matrices
and can be
defined as follows. Given $B\in \spn$, write the (unique) polar
decomposition of the matrix $e^{tB}$ as $e^{tB} = P(t)U(t)$. Here
$P(t)$ is symplectic (i.e. belongs to $Sp (2n, \R)$), symmetric
and positive and $U(t)$ is symplectic and complex-linear. The real
operator $U(t)$ can be identified with a {\it unitary} operator on
$\C^n$ and we denote by
$$det_{\C} U(t) \in S^1 =\{ z\in \C\ |\ |z|=1\}$$
the determinant of this unitary operator. Note that the families
$P(t), U(t)$ are continuous in $t$. Now set
$$\zeta_M (B) := \lim_{t\to +\infty} \frac{1}{t} \arg
det_{\C} U(t).$$

One can check that $\zeta_M$ is a continuous non-linear Lie
quasi-state -- in fact, it is a directional derivative of a unique
(up to a non-zero constant factor) homogeneous quasi-morphism on
$\tSp (2n,\R)$ (cf. Section~\ref{sec-Ad-inv-Lie-qs}).

For $n=1$ one can easily write an explicit formula for
$\zeta_M$. Namely,
$$sp\, (2,\R) = \Bigg\{ A=\left(
\begin{array}{cc}
  a & b \\
  c & -a \\
\end{array}
\right), a,b,c\in\R \Bigg\}$$ and
\begin{equation}
\label{eqn-Maslov-in-dim-2-def}
\zeta_M
(A) =
\begin{cases}
  \sqrt{|a^2 + bc|},\ {\rm if}\, a^2+bc <0, b<0, c>0, \\
  -\sqrt{|a^2 + bc|},\ {\rm if}\, a^2+bc <0, b>0, c<0, \\
  0,\ {\rm if}\, a^2+bc \geq 0. \\
\end{cases}
\end{equation}
As we see, $\zeta_M$ is continuous but not differentiable.

\subsection{Perspective and open questions}

Theorem~\ref{thm-main} raises the following general problem:
given a Lie algebra $\gog$, describe the space $\cQ (\gog)$
of continuous non-linear Lie quasi-states on $\gog$.

Besides the Lie algebras mentioned above, there are a few other
(finite-dimensional) cases where the answer is known: for
instance, the Heizenberg algebra (in this case $\cQ (\gog)=0$ --
see Section~\ref{subsec-heisenberg}) and the algebra
$\mathfrak{so} (3,\R)$ (in this case any two commuting elements
must be proportional to each other so any continuous
$\R$-homogeneous function is a continuous Lie quasi-state). A
partial result on Lie quasi-states on $\gogl$, needed for the
proof of Theorem~\ref{thm-main} can be found in Section
\ref{subsec-rank-1}. Otherwise, as far as the classical Lie
algebras are concerned, the answer is unknown already for
$\mathfrak{sl} (3,\R)$.

Further, in view of Theorem~\ref{thm-discont} it would be
interesting to relax the continuity assumption and to describe the
space of non-linear Lie quasi-states on $\gog$ that are {\it
bounded on a neighborhood of} $0$. We do not know the complete
answer even for the case of $\spn$. A possible interesting
modification of the question above would be to explore Lie
quasi-states on $\gog$ that are {\it positive on a certain cone}
in $\gog$. This is motivated by the original version of Gleason's
theorem and by the theory of symplectic quasi-states which have
this sort of property (see \cite{EP-qst}).

Another set of questions arises from the relation between Lie
quasi-states and homogeneous quasi-morphisms.

First, note that homogeneous quasi-morphisms appear as a part of a
certain remarkable cohomological theory on groups, called {\it
bounded cohomology} -- see e.g. \cite{Bav,Grigorchuk,Monod}. It
would be interesting to find a helpful cohomology theory for Lie
algebras incorporating Lie quasi-states.

Second, note that the directional derivative of a functional $\mu:
G\to \R$ on a Lie group $G$ is a Lie quasi-state on the Lie algebra
$\gog$ of $G$ provided
\begin{equation}
\label{eqn-generalized-qmm} \mu (xy) = \mu (x) + \mu (y) \ {\rm for\
all\ commuting}\ x,y\in G
\end{equation}
(and, in particular, $\mu (x^k) = k\mu(x)$ for any $k\in\Z$).   It
would be interesting to find whether a continuous $\mu$ satisfying
\eqref{eqn-generalized-qmm} always has to be a quasi-morphism and,
more generally, to describe the quotient of the space of all such
$\mu$ on a given $G$ by the space of continuous homogeneous
quasi-morphisms on $G$.

Third, given the Lie algebra $\gog$ of a (simply connected) Lie
group $G$, one can consider the following subsets of $\hcQ
(\gog)$:

$\hcQ_{qm} (\gog) := \{$the space of continuous Lie quasi-states
on $\gog$ coming from continuous homogeneous quasi-morphisms on $G
\}$;

$\hcQ_{Ad} (\gog) := \{$the set of continuous Lie quasi-states on
$\gog$ which are invariant under the adjoint action of $G  \}$.

Clearly,
$$ \hcQ_{qm} (\gog) \subset \hcQ_{Ad} (\gog) \subset \hcQ (\gog).$$
By Theorem~\ref{thm-main} these spaces coincide for $\spn$. It
would be interesting to explore these inclusions for other
algebras $\gog$. For instance, assume that $\gog$ is a compact
simple Lie algebra. In this case $\hcQ_{qm} (\gog)=0$ since any
continuous homogeneous quasi-morphism on a compact group has to be
zero. Further, we show in Section~\ref{sec-Ad-inv-Lie-qs} below
that $\hcQ_{Ad} (\gog)=0$. At the same time note that $\hcQ
(\gog)$ might sometimes be infinite-dimensional, for instance, if
$\gog = \gt{so} (3, \R)$.

Let us mention finally that the study of homogeneous
quasi-morphisms on groups is closely related to geometrical
structures and dynamics on spaces where these groups act -- see
e.g. \cite{Ghys}. It would be interesting to understand geometric
and/or dynamical meaning of non-linear Lie quasi-states, for
instance, those constructed in Section~\ref{sec-discont-Lie-qs}
below.

\section{Proof of the main theorem}\label{sec-proof-main}

Let $\zeta$ be a continuous Lie quasi-state on $\spn$, $n\geq 3$.
We want to show that $\zeta$ is a sum of $c\zeta_M$, $c\in\R$, and a linear
functional on $\spn$.

\subsection{Preliminaries}
\label{sec-prelim-main-thm}

For each $k\in\N$ denote by $M_k (\R)$ (respectively, $M_k (\C)$)
the spaces of real (respectively, complex) $k\times k$-matrices.

Let $\omega = \sum_{k=1}^n dp_k\wedge dq_k$ be the standard linear
symplectic form on the vector space $\R^{2n}$ with the coordinates
$p_1,\ldots,p_n, q_1,\ldots, q_n$ on it.

For each $A\in M_{2n} (\R)$ there exists a unique $A^\omega\in
M_{2n} (\R)$ such that $\omega (Ax,y) = \omega (x,Ay)$ for any
$x,y\in \R^{2n}$. We say that $A$ is {\it $\omega$-symmetric}, if
$A=A^\omega$, and {\it skew-symplectic}, if $A=-A^\omega$.
With this terminology, $\spn$ is the algebra of
skew-symplectic matrices $A\in M_{2n} (\R)$.

Given two vectors $\xi, \eta\in \R^{2n}$, define the following operators on $\R^{2n}$:
$$ T_{\xi,\eta} (x) = \omega (\xi, x)\eta,$$
$$Y_{\xi,\eta} (x) := T_{\xi,\xi} (x) + T_{\eta,\eta} (x)=\omega (\xi, x)\xi + \omega (\eta,x)\eta,$$
$$Z_{\xi,\eta} (x) := T_{\eta,\xi} (x) + T_{\xi,\eta} (x) = \omega (\eta, x)\xi + \omega (\xi,x)\eta.$$
One readily checks that
$$T_{\xi,\xi}, Y_{\xi,\eta}, Z_{\xi,\eta}\in \spn.$$
Note also that
$$Y_{\xi,\eta} = Y_{\eta,\xi},\ Z_{\xi,\eta} = Z_{\eta,\xi}$$
and $T_{\xi,\eta}, Z_{\xi,\eta}$ depend linearly on $\xi$ and
$\eta$. Finally, an easy computation shows that
$$\zeta_M (Y_{\xi,\eta}) = -|\omega(\xi,\eta)|,\ \zeta_M (Z_{\xi,\eta}) =0.$$

Consider $\C^{2n}= \R^{2n} \oplus i\R^{2n}$ as the complexification
of $\R^{2n}$. We write elements of $\C^{2n}$ as $v=a+ib$,
$a,b\in\R^{2n}$.

Denote by $(\cdot,\cdot )$ the standard Euclidean inner product on
$\R^{2n}$ and by $\langle \cdot, \cdot \rangle$ the standard
Hermitian inner product on $\C^{2n}$ so that
$$\langle a+ib, c+id\rangle = (a,c) + (b,d) + i (b,c) - i(a,d).$$

We say that a function $Q: \C^{2n}\to\R$ is a {\it real Hermitian
quadratic form}  if $Q$ is a real quadratic form satisfying $Q
(\lambda v) = |\lambda|^2 Q (v)$ for any $\lambda\in\C$, or,
equivalently, if $Q (v) = \cH(v,v)$ for some Hermitian form $\cH:
\C^{2n}\times\C^{2n}\to\C$. By definition, a Hermitian form $\cH$ is
given by $\cH (v,w) = \langle Hv, w\rangle$ for some Hermitian
$2n\times 2n$-matrix $H$, which, in turn, can be always written as
$H= A+iB$, where $A, B\in M_{2n} (\R)$, $A=A^T$, $B= - B^T$. The
corresponding real Hermitian quadratic form $Q$ can then be written
as
\begin{equation}
\label{eqn-Q-inner-product}
Q (a+ib) = (Aa,a) + (Ab,b) + 2(Ba,b).
\end{equation}
On the other hand, since $\omega$ is non-degenerate, any real
bilinear form on $\R^{2n}$ can be uniquely represented as $\omega
(C\cdot, \cdot)$ for some $C\in M_{2n} (\R)$, and moreover, such a
bilinear form is symmetric (respectively, anti-symmetric) if and
only if $C=-C^\omega$ (respectively, $C=C^\omega$). Together with
\eqref{eqn-Q-inner-product} this yields that any real Hermitian
quadratic form $Q$ on $\C^{2n}$ can be written as
\begin{equation}
\label{eqn-Q-omega} Q (a+ib) = \omega (Ca,a) +\omega (Cb,b) +
\omega (Da,b),
\end{equation}
where
$$C,D\in M_{2n} (\R),\, C=-C^\omega,\, D=D^\omega.$$

\subsection{Reduction to the computation of $\zeta$ on $Y_{\xi,\eta},
Z_{\xi,\eta}$}

First, we will show the following

\begin{prop}
Any continuous Lie quasi-state $\zeta$ on $\spn$ is completely
determined by its values on elements of the form $Y_{\xi,\eta},
Z_{\xi,\eta}$.
\end{prop}

\begin{proof}


Since semi-simple (i.e. diagonalizable over
$\C$) elements are dense in $\spn$ and $\zeta$ is continuous, it is enough to
show that for any semi-simple $A\in\spn$ the computation of
$\zeta(A)$ can be reduced to the computation of $\zeta$ on some $Y_{\xi,\eta},
Z_{\xi,\eta}$.

Recall that a Darboux basis $e_1,\ldots, e_n, f_1,\ldots, f_n$ on
$\R^{2n}$ is a basis which satisfies
$$\omega(e_i,e_k) = \omega(f_i,f_k) = 0,\;\; \omega(e_i,f_k) =
\delta_{ik}\;,$$ for all $i,k =1,...,n$. Here $\delta_{ik} =0$ for
$i\neq k$ and $\delta_{ii}=1$.

By Williamson's results on the normal forms of skew-symplectic
matrices with respect to the adjoint action of the symplectic
group (see \cite{Williamson}, cf. \cite{Arnold-math-methods},
Appendix 6, and \cite{Hormander}), any semi-simple $A\in\spn$ can
be represented in an appropriate basis as a block-diagonal matrix,
where the blocks correspond to symplectic subspaces of $\R^{2n}$
spanned by vectors $e_k,\ldots, e_{k+l}, f_k,\ldots , f_{k+l}$
from a Darboux basis. Each such block, written in the basis
$e_k,\ldots, e_{k+l}$, $f_k,\ldots , f_{k+l}$, belongs to one of
the following three types:

\begin{itemize}

\item[(i)]
$\left(%
\begin{array}{cc}
  -a & 0 \\
  0 & a \\
\end{array}%
\right)$,
corresponding to a $2$-dimensional symplectic subspace of $\R^{2n}$ spanned by some $e_k, f_k$.

\item[(ii)]
$\left(%
\begin{array}{cc}
  0 & b \\
  -b & 0 \\
\end{array}%
\right)$,
corresponding to a $2$-dimensional symplectic subspace of $\R^{2n}$ spanned by some $e_k, f_k$.

\item[(iii)]
$\left(%
\begin{array}{cccc}
  -a & b & 0 & 0 \\
  -b & -a & 0 & 0 \\
  0 & 0 & a & b \\
  0 & 0 & -b & a \\
\end{array}%
\right)$,
corresponding to a $4$-dimensional symplectic subspace of $\R^{2n}$ spanned by some $e_k, e_{k+1}, f_k, f_{k+1}$.

\end{itemize}

Note that
\begin{itemize}

\item[I.] A block-diagonal $2n\times 2n$-matrix having only one
block which is the $2\times 2$-block (i) can be represented as $a
Z_{e_k, f_k}$.

\item[II.] A block-diagonal $2n\times 2n$-matrix having only one
block which is the $2\times 2$-block (ii) can be represented as $b
Y_{e_k, f_k}$.

\item[III.] A block-diagonal $2n\times 2n$-matrix $A$ having only one
block which is the $4\times 4$-block (iii) can be represented as
$$A=a Z_{e_k, f_k} + a Z_{e_{k+1}, f_{k+1}} -
b Y_{\frac{e_{k+1}-f_k}{\sqrt{2}},\frac{e_k+f_{k+1}}{\sqrt{2}}} +
b Y_{\frac{e_k-f_{k+1}}{\sqrt{2}}, \frac{e_{k+1}+f_k}{\sqrt{2}}},$$
where
\begin{equation}
\label{eqn-commut-Y-Z-1} \bigg[a Z_{e_k, f_k} + a Z_{e_{k+1},
f_{k+1}}, -b Y_{\frac{e_{k+1}-f_k}{\sqrt{2}},
\frac{e_k+f_{k+1}}{\sqrt{2}}} + b Y_{\frac{e_k-f_{k+1}}{\sqrt{2}},
\frac{e_{k+1}+f_k}{\sqrt{2}}}\bigg] =0,
\end{equation}
\begin{equation}
\label{eqn-commut-Y-Z-2}
[Z_{e_k, f_k}, Z_{e_{k+1}, f_{k+1}}]=0,
\end{equation}
\begin{equation}
\label{eqn-commut-Y-Z-3}
\bigg[Y_{\frac{e_{k+1}-f_k}{\sqrt{2}},\frac{e_k+f_{k+1}}{\sqrt{2}}},
Y_{\frac{e_k-f_{k+1}}{\sqrt{2}}, \frac{e_{k+1}+f_k}{\sqrt{2}}}\bigg]=0.
\end{equation}

\end{itemize}

Thus every semi-simple $A\in\spn$ can be represented as a sum of
{\it pairwise commuting} block-diagonal matrices, each of which
has one of the forms I, II, III, and thus the computation of
$\zeta (A)$ reduces to the computation of $\zeta$ on the
block-diagonal matrices I, II, III. The latter in turn reduces to
the computation of $\zeta$ on some $Y_{\xi, \eta}$, $Z_{\xi,
\eta}$ -- for I and II this is trivial and for III this
follows from the commutation relations \eqref{eqn-commut-Y-Z-1},
\eqref{eqn-commut-Y-Z-2}, \eqref{eqn-commut-Y-Z-3}.

\end{proof}

Given $\xi, \eta\in \R^{2n}$, denote
\begin{equation}
\label{eqn-def-F} F(\xi,\eta):=\zeta (Y_{\xi,\eta}),
\end{equation}
\begin{equation}
\label{eqn-def-G} G(\xi,\eta) := \zeta (Z_{\xi,\eta}).
\end{equation}

Let $\C^{2n}=\R^{2n} + i\R^{2n}$ be the complexification of
$\R^{2n}$. Put
$$\cW:=\{ x + iy\in \C^{2n}, \omega (x,y)> 0\}\;.$$
Consider the function $F(\xi,\eta)$ as a function on $\C^{2n}$:
$F(\xi,\eta):=F(\xi+i\eta)$. We will prove the following key
proposition:

\medskip
\noindent
\begin{prop}
\label{prop-F-G-quadratic} $\;$

\begin{itemize} \item[{(i)}]There exists a real Hermitian quadratic
form on $\C^{2n}$ which coincides with $F$ on $\cW$.
\item[{(ii)}] The function $G$ is a symmetric bilinear form on
$\R^{2n}$.\end{itemize}

\end{prop}

The proof of the proposition uses Gleason's theorem, and hence the
assumption $n\geq 3$. Postponing the proof,
let us show how the proposition implies the main theorem.

\medskip
\noindent {\bf Deducing Theorem~\ref{thm-main} from Proposition
\ref{prop-F-G-quadratic}:} By \eqref{eqn-Q-omega} we can write the
real Hermitian quadratic form $F$ as
\begin{equation}
\label{eqn-F-omega-C-D} F(\xi+i\eta)= \omega (C\xi, \xi) + \omega
(C\eta, \eta) + \omega (D\xi,\eta) \;\forall \; \xi+i\eta \in \cW
\end{equation}
for some $C,D \in M_{2n} (\R)$, $C=-C^\omega, D=D^\omega$. Put
$$\cW^-:=\{ x + iy\in \C^{2n}, \omega (x,y)< 0\}\;.$$
Observe that $Y_{\xi,\eta} = Y_{\xi,-\eta}$ and hence, since $\zeta$
is $\R$-homogeneous,
$$F(\xi+i\eta) = F(\xi- i\eta)\; \forall \xi,\eta \in \R^{2n}\;.$$
Since $\xi+i\eta \in \cW$ whenever $\xi-i\eta \in \cW^-$ we get that
$$F(\xi+i\eta)= \omega (C\xi, \xi) + \omega (C\eta, \eta) - \omega
(D\xi,\eta) \;\forall \xi+i\eta \in \cW^-\;.$$

Next, we claim that $D = c\id$ for some $c \in \R$. Indeed,
otherwise there exist $\xi,\eta$ so that $\omega (D\xi,\eta) \neq
0$ and $\omega(\xi,\eta)=0$. Hence $F$ is discontinuous at $\xi
+i\eta$, so we get a contradiction.

Since $\cW \cup \cW^-$ is dense in $\C^{2n}$ and $F$ is
continuous, we conclude that
\begin{equation}\label{eq-explicit-F}
F(\xi+i\eta) = \omega (C\xi, \xi) + \omega (C\eta, \eta) + c|\omega
(\xi,\eta)| \;\;\forall \; \xi+i\eta \in \C^{2n}\;.\end{equation}

Noting that $Y_{\xi,\xi} = Z_{\xi,\xi}$ (and hence $F(\xi+i\xi) =
\zeta (Y_{\xi,\xi}) = \zeta (Z_{\xi,\xi}) = G(\xi, \xi)$) we get
$$F(\xi+i\xi)= G(\xi,\xi) = 2\omega (C\xi, \xi),$$
and therefore
 \begin{equation}\label{eq-explicit-G} G(\xi,\eta) =
2\omega (C\xi, \eta)\;.\end{equation}

Define a linear functional $\alpha: \spn \to\R$ by $\alpha (A):=
-tr (CA)$. We claim that $\zeta$ equals to $\alpha - c\zeta_M$ on
each $Y_{\xi,\eta}$, $Z_{\xi,\eta}$ -- as it was explained in the
beginning of this section, this claim would imply the theorem.

In order to prove the claim, we observe that
$$\alpha (Y_{\xi,\eta}) = -tr (CY_{\xi,\eta}) = -tr (CT_{\xi,\xi}) - tr (CT_{\eta,\eta}) =$$
$$= -tr T_{\xi, C\xi} - tr T_{\eta, C\eta} = \omega (C\xi, \xi) + \omega (C\eta, \eta),$$

$$\alpha (Z_{\xi,\eta}) = -tr (CZ_{\xi,\eta}) = -tr (CT_{\xi,\eta}) - tr (CT_{\eta,\xi}) =$$
$$ = -tr T_{\xi, C\eta} - tr T_{\eta, C\xi} = \omega (C\eta, \xi) + \omega (C\xi, \eta) = 2\omega (C\xi, \eta).$$
Thus, by \eqref{eq-explicit-F} and \eqref{eq-explicit-G}
$$\zeta (Y_{\xi,\eta}) = \alpha (Y_{\xi,\eta}) + c|\omega (\xi,\eta)|,$$
$$\zeta (Z_{\xi,\eta}) = \alpha (Z_{\xi,\eta})\;.$$
The claim follows, since
$$\zeta_M (Y_{\xi,\eta}) = -|\omega(\xi,\eta)|,\ \zeta_M (Z_{\xi,\eta}) = 0.$$
This finishes the proof of the theorem modulo
Proposition~\ref{prop-F-G-quadratic}.\Qed

\bigskip
\noindent{\bf Outline of the proof of
Proposition~\ref{prop-F-G-quadratic}(i):} We reduce the problem to
the case when $F$ is smooth on $\cW$ (see Section
\ref{subsec-smooth}). Further, for any $\omega$-compatible almost
complex structure $J$ on $\R^{2n}$ the space $L_J:={\xi +iJ\xi}$
is Lagrangian with respect to the complexification $\omega^{\C}$
of $\omega$. We observe that the transformation $Y_{\xi,J\xi}$
lies in the $J$-unitary subalgebra of $\spn$, and hence Gleason's
theorem (complex version) yields that the restriction of $F$ to
each $L_J$ is a Hermitian quadratic form, say, $q_J$. Note that
$q_J=q_I$ on $L_J \cap L_I$  for every two compatible almost
complex structures $J,I$. This yields a restriction on the
differential ${\partial q_J}/{\partial J}$ which can be translated
into a system of first order PDE's. The analysis of the system
(see Section~\ref{sec-functions-whose-restr-are-qforms})
eventually yields that $F$ on $\cW$ is quadratic (see
Section~\ref{subsec-proof-i}).

\medskip
\noindent{\bf Outline of the proof of
Proposition~\ref{prop-F-G-quadratic}(ii):} The proof of (ii) is a
bit trickier. First, we show that for $n\geq 3$ any continuous Lie
quasi-state on $\gogl$, when restricted to rank 1 operators
$B_{\xi,\eta}x= (x,\eta)\xi$, is given by $\zeta(B_{\xi,\eta})= tr
NB_{\xi,\eta}$ for some fixed matrix $N$. Interestingly enough, the
proof of this statement is very similar to the proof of
Proposition~\ref{prop-F-G-quadratic}(i) outlined above, but {\it
over the field} $\R$, see Section~\ref{subsec-rank-1}.

This result readily yields that for a fixed $\xi$, the restriction
of $G(\xi,\cdot)$ to every Lagrangian subspace of $\R^{2n}$ is
linear, see Section~\ref{subsec-proof-ii}.

Finally, we make a detour to the Heisenberg Lie algebra: we show
that on this algebra every (not necessarily continuous)
quasi-state is linear, see Section~\ref{subsec-heisenberg}. As an
immediate consequence we get that $G(\xi,\eta)$ is linear in the
variable $\eta$. Since $G$ is symmetric in $\xi$ and $\eta$, this
completes the proof.

\bigskip

The rest of Section~\ref{sec-proof-main} contains a proof of
Proposition~\ref{prop-F-G-quadratic}. The plan described above
serves as a rough guideline only, the specific details are often
formulated in a different language and appear in a different
order.

\subsection{Smoothening Lie quasi-states}\label{subsec-smooth}
Eventually, we wish to reduce Proposition~\ref{prop-F-G-quadratic}
to the case when the functions $F$ and $G$ are smooth on $\cW\cup
\cW^-$. For that purpose we show  that a continuous Lie
quasi-state on any Lie algebra can be suitably approximated by Lie
quasi-states which are smooth along the orbits of the adjoint
action of the Lie group. We thank Semyon Alesker for explaining
this to us.

Let $\gog$ be the Lie algebra of a (connected) Lie group $G$. Fix a
norm $\|\cdot\|$ on $\gog$, set
$$S:= \{ x\in\gog, \| x\|=1\}$$
and define a metric $d$ on
$\hcQ (\gog)$ by
$$d (\zeta_1, \zeta_2) := \sup_{x\in S} |\zeta_1 (x) - \zeta_2
(x)|.$$
Limits with respect to $d$ will be called $d$-limits.

We will say that a function $\zeta: \gog\to\R$ is {\it
orbit-smooth} if the restriction of $\zeta$ on any orbit of the
adjoint action of $G$ on $\gog$ is smooth, or, in other words, if
the function $g\mapsto \zeta(gx)$ is a smooth function on $G$ for
every $x\in \gog$. (Here $gx$ denotes the adjoint action of $g$ on
$x$).

\begin{prop}
\label{prop-smoothening} Every $\zeta\in \hcQ (\gog)$ is a $d$-limit
of orbit-smooth continuous Lie
quasi-states.

\end{prop}

\begin{proof}
Let $\mu$ be a right-invariant smooth measure on $G$. We will
measure diameters of sets in $G$ with respect to some fixed distance
function defined by a Riemannian metric on $G$.

 Let $\varphi_i$, $i\in\N$, be a delta-like sequence of
$C^\infty$-smooth functions on $G$. In particular, assume that
\begin{itemize}
\item{} $\varphi_i\geq 0$ for all $i$;

\item{} $\supp \varphi_i\subset U_i$, $i\in\N$, for some open neighborhoods $U_i$ of $\id\in G$ such that
$\diam U_i\to 0$ as $i\to +\infty$;

\item{} $\int_G \varphi_i d\mu = 1$ for all $i$.

\end{itemize}

Given $i\in\N$, put $$\zeta_i (x):= \int_G \zeta (hx)\varphi_i (h)
d\mu (h).$$ Note that since $\zeta$ is continuous, $\zeta_i:
\gog\to\R$ is continuous as well. Moreover, for any $x_0\in \gog$
and $g\in G$ we have
$$\zeta_i (gx_0) = \int_G \zeta (hgx_0) \varphi_i (h) d\mu (h) =
\int_G \zeta (h'x_0)\varphi_i (h'g^{-1})d\mu (h'g^{-1}),$$
where $h'=hg$. Since $\mu$ is right-invariant, the latter integral
is equal to
$$\int_G
\zeta (h' x_0)\varphi_i (h'g^{-1})d\mu (h').$$ Obviously, this
integral depends smoothly on $g$. Thus for any $i\in\N$ the
functional $\zeta_i$ is orbit-smooth.

Let us estimate $d (\zeta_i, \zeta)$. Note that
$$|\zeta_i (x) - \zeta (x)| = \int_G |\zeta (hx) - \zeta (x)|
\varphi_i (h) d\mu (h) =$$
$$= \int_{U_i} |\zeta (hx) - \zeta (x)| \varphi_i (h) d\mu (h)
\leq \max_{h\in U_i} |\zeta (hx) - \zeta (x)|.$$
Since $\diam U_i \to 0$ and $\zeta$ is uniformly continuous on $S$,  there exist $\delta_1 (i)$, $0< \delta_1 (i) <
1$,
and $\delta_2 (i)>0$, such that $\lim_{i\to +\infty} \delta_1 (i) =
0$, $\lim_{i\to +\infty} \delta_2 (i) = 0$, and such that
for any $x\in S$ and any $h\in U_i$
$$1- \delta_1 (i) \leq \| hx\| \leq 1+ \delta_1 (i),$$
$$\bigg| \zeta \bigg(\frac{hx}{\| hx\|}\bigg) - \zeta (x)\bigg|\leq \delta_2 (i).$$
Then, by the homogeneity of $\zeta$, we get for any $x\in S$ and
any $h\in U_i$
$$| \zeta (hx) - \zeta (x)| = \bigg| \| hx\| \zeta \bigg(\frac{hx}{\|hx\|}\bigg) - \zeta
(x)\bigg|\leq$$
$$\leq \bigg| \|hx\| \zeta \bigg(\frac{hx}{\|hx\|}\bigg) - \|hx\| \zeta (x) \bigg| + \bigg|\|hx\|
\zeta (x) - \zeta (x)\bigg|\leq $$
$$\leq \|hx\|\delta_2 (i) + |\zeta(x)|\delta_1
(i)\leq \bigg(1+\delta_1 (i)\bigg)\delta_2 (i) + \max_S |\zeta| \cdot \delta_1
(i).$$
Thus
$$d (\zeta_i, \zeta) = \max_{x\in S} | \zeta (hx) - \zeta (x)| \leq
\bigg(1+\delta_1 (i)\bigg)\delta_2 (i) + \max_S |\zeta|\cdot\delta_1
(i),$$
and therefore
$\lim_{i\to +\infty} d (\zeta_i, \zeta) = 0$.
\end{proof}

\subsection{Functions whose restrictions on Lagrangian subspaces are
qua\-dratic forms} \label{sec-functions-whose-restr-are-qforms}

Recall that
$$\cW =\{ x + iy\in \C^{2n}, \omega (x,y)> 0\}\;.$$
Denote by $\cJ$  the space of all $J\in M_{2n} (\R)$ such that
$J^2=-\id$ and $(\cdot,\cdot)_J :=\omega (\cdot, J\cdot)$ is a
$J$-invariant inner product on $\R^{2n}$. (Such a $J$ is called a
{\it a complex structure on $\R^{2n}$ compatible with $\omega$}).

Given $J\in \cJ$, define a complex vector subspace $L_J$ of
$\C^{2n}$ by $$L_J:=  \{ x +iJx\in \C^{2n},\ x\in \R^{2n}\}.$$

Consider a complex-valued symplectic form $\omega^\C$ on $\C^{2n} =
\R^{2n} \oplus i\R^{2n}$  which is the complexification of $\omega$:
$$\omega^\C (a+ib, c+id) = \omega (a,c) -\omega (b,d) + i \big(\omega
(b,c) +\omega (a,d)\big).$$ Note that each $L_J$, $J\in\cJ$, is
Lagrangian with respect to $\omega^\C$ and, more generally, an
$\omega^\C$-Lagrangian complex vector subspace $L\subset \C^{2n}$
has the form $L=L_J$ for some $J\in\cJ$  if and only if $L\setminus
0\subset \cW$. The set of the subspaces $L_J$, $J\in\cJ$, is open in
the set of $\omega^\C$-Lagrangian complex vector subspaces of
$\C^{2n}$.

\begin{prop}
\label{prop-complex-qforms-smooth-restrictions-to-LJ}

Let $F: \cW \to\R$ be a $C^3$-smooth  function. Assume that for
any $\omega^\C$-Lagrangian complex subspace $L_J\subset \cW$,
$J\in\cJ$, the restriction of $F$ to $L_J\cap \cW$  is a real
Hermitian quadratic form. Then $F$ is the restriction to $\cW$ of
a real Hermitian quadratic form on $\C^{2n}$.

\end{prop}

The proof will be based on the following proposition. Denote the
space of symmetric complex $n\times n$-matrices by $\cS_n (\C)$.
Let $(z,w)$ be  complex linear coordinates on the vector space
$\C^{2n}= \C^n\times\C^n$, where $z=(z_1,\ldots z_n),
w=(w_1,\ldots, w_n)$. For an open connected neighborhood $\cV$ of
$0$ in $\cS_n (\C)$ put
$$\cC_{\cV} := \{(z,Az) \in \C^{2n}\;:\; z \neq 0, A \in \cV\}\;.$$
One readily checks that the set $\cC_{\cV}$ is open and invariant
under multiplication by non-zero scalars from $\C$.

\begin{prop}
\label{prop-complex-qforms-restrictions-to-LA}  Let $\cV$ be an open
connected neighborhood of $0$ in  $\cS_n (\C)$, $n \geq 3$. Let $F:
\cC_{\cV} \to \R$ be a $C^3$-smooth function so that
\begin{equation}
\label{eqn-F-homogen} F (\lambda v) = |\lambda|^2 F (v)\
\forall\lambda\in\C \setminus\{0\},\ \forall v\in  \cC_{\cV} \;.
\end{equation}
Assume that the restriction of $F$ to any vector subspace
$$L_A:= \{ w= Az\}\subset \C^{2n},\ A\in\cV,$$
is a real Hermitian quadratic form. Then the function $F$ is the
restriction to  $\cC_{\cV}$ of some real Hermitian quadratic form on
$\C^{2n}$.
\end{prop}

As one can easily see, Proposition~\ref{prop-complex-qforms-smooth-restrictions-to-LJ}
and
Proposition~\ref{prop-complex-qforms-restrictions-to-LA} fail for
$n=1$. Nevertheless they do hold for $n=2$ though in this case one needs to
modify the proof of Proposition~\ref{prop-complex-qforms-restrictions-to-LA}
slightly -- see Remark~\ref{rem-F-A-n-equal-2} below.

\medskip
\noindent {\bf Deducing
Proposition~\ref{prop-complex-qforms-smooth-restrictions-to-LJ}
from Proposition~\ref{prop-complex-qforms-restrictions-to-LA}:}

Pick any $L:=L_{J_0}$, $J_0\in\cJ$. Using the linear Darboux
theorem for complex symplectic forms choose complex coordinates
$z=(z_1,\ldots z_n), w=(w_1,\ldots, w_n)$ on the vector space
$\C^{2n}$ so that $\omega^\C = dz\wedge dw$ and $L = \{ w=0\}$.
Fix a sufficiently small open connected neighborhood $\cV$ of zero
in $\cS_n (\C)$. Then $\cC_\cV\subset \cW$ (since $L\setminus
0\subset \cW$) and $F$ is $C^3$-smooth on $\cC_\cV$ (since it is
$C^3$-smooth on $\cW$). Moreover, any $L_A$, $A\in \cV$, has the
form $L_A=L_J$ for some $J\in\cJ$. Therefore, by
Proposition~\ref{prop-complex-qforms-restrictions-to-LA}, $F$
coincides on $\cC_{\cV}$ with the restriction of some real
Hermitian quadratic form defined on $\C^{2n}$. Now letting $J_0$
vary inside $\cJ$ we see that $\cW$ can be covered by open cones,
invariant under the multiplication by non-zero complex scalars, on
each of which $F$ coincides with the restriction of a real
Hermitian quadratic form. Since $F$ is $C^3$-smooth on $\cW$ and
$\cW$ is path-connected, this yields that $F$ coincides on the
whole $\cW$ with the restriction of some real Hermitian quadratic
form defined on $\C^{2n}$. \Qed

\bigskip \noindent {\bf Proof of
Proposition~\ref{prop-complex-qforms-restrictions-to-LA}:} Represent
$\C^{2n}$ as $\C^{2n} = \C^n\times \C^n$, where $z$ and $w$ are
coordinates along, respectively, the first and the second factors.
Accordingly, we will write the vectors in $\C^{2n}$ in the form
$z\oplus w$. Given $A\in \cV$, write
$$F(z,Az) = \langle\langle H(A) z, z\rangle\rangle,$$
where $H(A)$ is a Hermitian $n\times n$-matrix and
$\langle\langle\cdot,\cdot\rangle\rangle$ is the standard
Hermitian inner product on $\C^n$. Since $F$ is $C^3$-smooth, the
matrix $H(A)$ depends $C^3$-smoothly on  $A\in\cV$.

We want to show that there exists a Hermitian $2n\times 2n$-matrix
$\cH$ such that for any $A\in\cV$ and any $z$
\begin{equation}
\label{eqn-F-cH-A} F(z, Az) = \langle \cH (z\oplus Az), z\oplus
Az\rangle,
\end{equation}
$\langle\cdot,\cdot\rangle$ is the standard Hermitian inner
product on $\C^{2n}$. Write $\cH$ as a matrix with four $n\times
n$-blocks:
$$\cH = \left(%
\begin{array}{cc}
  P & Q \\
  \bar{Q}^T & R \\
\end{array}%
\right).$$
Note that $P$ and $R$ are Hermitian $n\times n$-matrices.
Rewriting \eqref{eqn-F-cH-A} we see that
we need to show that
$$H(A) = P + QA + \bar{A} \bar{Q}^T + \bar{A}RA,$$
or, in terms of matrix coefficients,
\begin{equation}
\label{eqn-F-cH-A-matrix-coord}
H_{ij} (A) = P_{ij} + \sum_\alpha Q_{i\alpha} A_{\alpha j} +
\sum_\alpha \bar{A}_{i\alpha}\bar{Q}_{j\alpha} +
\sum_{\alpha,\beta} \bar{A}_{i\alpha} R_{\alpha\beta} A_{\beta
j},\ \forall i,j=1,\ldots,n.
\end{equation}

\begin{rem}
{\rm
The coordinates on the space of symmetric matrices are $A_{ij}$
with $i\leq j$. Nevertheless, we shall also use the coordinates
$A_{ij}$ with $i>j$ identifying them with $A_{ji}$: $A_{ij}=
A_{ji}$.

}
\end{rem}

For any $1\leq s < t\leq n$ put
$$u_{ts} (z) = u_{st} (z) = z_s z_t,\ u_{ss} (z) = - z_t^2,\
u_{tt} (z) = -z_s^2.$$
Define a matrix $V^{s,t} (z) = (V^{s,t}_{ij} (z))_{i,j=1,\ldots,n}$
by
$$V^{s,t}_{ij} (z) =
\begin{cases}
u_{ij} (z),\ {\rm if}\ i,j\in \{ s,t\},\\
0,\ {\rm otherwise}.
\end{cases}
$$
Note that for any $z$
$$V^{s,t} (z) z = 0.$$
Thus for any $s,t$ the expression
$$\langle H(A +\epsilon V^{s,t} (z)) z, z\rangle = F (z, Az +\epsilon V^{s,t}
(z)z) = F(z,Az)$$
does not depend on $\epsilon$, where $\epsilon$ is a complex parameter.
Differentiating by $\epsilon$ at $\epsilon=0$ we get the following
system of equations for every $i,j,s,t$, $s<t$, and any $A\in \cV$:
$$\sum_{i,j} \bigg( \frac{\partial H_{ij}}{\partial A_{st}} (A) u_{st} (z) z_j \bar{z}_i
+ \frac{\partial H_{ij}}{\partial A_{ss}} (A) u_{ss} (z) z_j
\bar{z}_i+
\frac{\partial H_{ij}}{\partial A_{tt}} (A) u_{tt} (z) z_j
\bar{z}_i+$$
$$
+
\frac{\partial H_{ij}}{\partial \bar{A}_{st}} (A) \bar{u}_{st} (z) z_j \bar{z}_i
+
\frac{\partial H_{ij}}{\partial \bar{A}_{ss}} (A) \bar{u}_{ss} (z) z_j \bar{z}_i
+
\frac{\partial H_{ij}}{\partial \bar{A}_{tt}} (A) \bar{u}_{tt} (z) z_j
\bar{z}_i \bigg) = 0
$$
for any $z, \bar{z}$.
This is a polynomial in $z,\bar{z}$ and, since it vanishes on an
open set, all its coefficients have to vanish. These coefficients
can be found by collecting similar terms in the last equation.
A straightforward analysis of these terms yields that for any
$i,j,s,t=1,\ldots,n$, $s<t$, the partial derivatives of $H_{ij}$ at any
$A\in\cV$ satisfy the following equations:
\begin{equation}\label{eq-part-1}
\frac{\partial H_{it}}{\partial A_{st}} =
\frac{\partial H_{is}}{\partial A_{ss}},
\end{equation}
\begin{equation}\label{eq-part-2}
\frac{\partial H_{is}}{\partial A_{st}} =
\frac{\partial H_{it}}{\partial A_{tt}},
\end{equation}
\begin{equation}\label{eq-part-3}
\frac{\partial H_{tj}}{\partial \bar{A}_{st}} = \frac{\partial
H_{sj}}{\partial \bar{A}_{ss}},
\end{equation}
\begin{equation}\label{eq-part-4}
\frac{\partial H_{sj}}{\partial \bar{A}_{st}} = \frac{\partial
H_{tj}}{\partial \bar{A}_{tt}}.
\end{equation}
Furthermore,
\begin{equation}\label{eq-part-5}
\frac{\partial H_{ij}}{\partial A_{\alpha\beta}} = 0\ {\rm if}\
j\notin\{ \alpha,\beta\},
\end{equation}
\begin{equation}\label{eq-part-6}
\frac{\partial H_{ij}}{\partial \bar{A}_{\alpha\beta}} = 0\ {\rm
if}\ i\notin\{ \alpha,\beta\}.\end{equation}
Note that equations
\eqref{eq-part-1} and \eqref{eq-part-2} can be summarized as
$$
\frac{\partial H_{ij}}{\partial A_{lj}} = \frac{\partial
H_{ir}}{\partial A_{lr}}\;\; \forall i,j,l,r\;.
$$
Differentiating this equation by $A_{kj}$ and using
\eqref{eq-part-5} one gets that
$$\frac{\partial^2 H_{ij}}{\partial A_{kj}\partial A_{lj}} = \frac{\partial^2 H_{ir}}{\partial A_{kj}\partial
A_{lr}}= 0,$$ provided $r\notin\{ k,j\}$. But such an $r$ always
exists since $n\geq 3$. Thus
\begin{equation}\label{eq-part-500}
\frac{\partial^2 H_{ij}}{\partial A_{kj}\partial A_{lj}} = 0\
\forall i,j,k,l\;.\end{equation} Similarly, using equations
\eqref{eq-part-3},\eqref{eq-part-4} and \eqref{eq-part-6} one gets
that\begin{equation}\label{eq-part-501}
\frac{\partial^2 H_{ij}}{\partial \bar{A}_{ik}\partial
\bar{A}_{il}} = 0\
\forall i,j,k,l\;.\end{equation}

Observe now that for all $i,j$  all the third derivatives of
$H_{ij}$ with respect to $A_{st}, \bar{A}_{st}$ vanish. Indeed, in
any third derivative either $A$-variables or $\bar{A}$-variables
appear at least twice, and the result follows from the vanishing of
the corresponding lower order derivatives, see formulas \eqref{eq-part-5}--\eqref{eq-part-501}.
Hence, each $H_{ij}$, as a
function of the variables $A_{st}, \bar{A}_{st}$, is a
(non-homogeneous) quadratic polynomial. The equations above on the
first and second partial derivatives of $H_{ij}$ allow to recover
the coefficients of this quadratic polynomial and check that this
polynomial indeed has the form \eqref{eqn-F-cH-A-matrix-coord}. \Qed

\begin{rem}
\label{rem-F-A-n-equal-2} {\rm In the case $n=2$, a slightly more
fine analysis of equations \eqref{eq-part-1}-\eqref{eq-part-6}
yields the same result: $H_{ij}$'s are non-homogeneous quadratic
polynomials of variables $A_{st}, \bar{A}_{st}$ of the form
\eqref{eqn-F-cH-A-matrix-coord}. }
\end{rem}

A completely similar argument yields the following proposition,
which we will need later and which is an analogue of
Proposition~\ref{prop-complex-qforms-restrictions-to-LA} for
functions on a real vector space.

Let $(\cdot, \cdot)$ be the Euclidean inner product on $\R^n$.
Set
$$\cU :=\{ x\times y\in \R^{2n}= \R^n\times \R^n\ | (x,y) > 0\}.$$
Denote by $\cS^+_n (\R)$ the space of symmetric real
positive-definite $n\times n$-matrices.

\begin{prop}
\label{prop-real-qforms-restrictions-LA}

Let $F: \R^n\times \R^n\to\R$, $n\geq 2$, be a continuous function
which is $C^3$-smooth on $\cU$. Assume that the restriction of $F$
on any vector subspace $$L_A:= \{ y= Ax\}\subset \R^n\times \R^n,\
A\in\cS^+_n (\R),$$ is a quadratic form. Then there exists a
quadratic form $Q$ on $\R^{2n}$ which coincides with $F$ on $\cU$.
\end{prop}


\medskip
\noindent\begin{rem}\label{rem-quadr}{\rm Denote by $\cL$ the
Lagrangian Grassmannian of the symplectic vector space $\R^{2n}$.
For a connected open subset $U \subset \cL$ consider the set
$\cC:=\bigcup_{L \in U} L \setminus \{0\}$. Let $F:\cC \to \R$ be
a continuous function whose restriction to every $L \in U$ is a
quadratic form. {\it Is it true that $F$ is the restriction of
some quadratic form defined on $\R^{2n}$ to $\cC$?} The analogue
over $\R$ of Proposition
\ref{prop-complex-qforms-restrictions-to-LA} above yields  the
affirmative answer provided $n \geq 3$ and $F$ is $C^3$-smooth. We
already mentioned that a small modification of our argument
settles the $n=2$ case, and there is a strong evidence that
$C^2$-smoothness suffices as well. It would be interesting to
understand what is precisely the minimal regularity assumption on
$F$ for which the question admits the positive answer.

In case when $U =\cL$, that is when $F$ is defined on the whole
$\R^{2n}$, the answer is affirmative: one uses the smoothening in
the spirit of Section~\ref{subsec-smooth} in order to reduce the
problem to the case when $F$ is not only continuous on $\R^{2n}$
but also smooth on $\R^{2n}\setminus 0$. Interestingly enough,
even when $U=\cL$ there exist discontinuous $F$ whose restrictions
to every Lagrangian subspace are quadratic forms. These examples
were constructed by Gleason in \cite{Gleason-2}.

Let us mention that this circle of problems extends verbatim into
the complex setting.}\end{rem}

\subsection{Proof of
Proposition~\ref{prop-F-G-quadratic}(i).}\label{subsec-proof-i}

Note that the group $Sp (2n,\R)$ acts transitively on the set of
pairs $\xi\times\eta\in \R^{2n}\times \R^{2n}$ such that $\omega
(\xi,\eta) = 1$. Hence, the set $\{ Y_{\xi, \eta}\ | \
\omega(\xi,\eta)=1\}$ is an orbit of the adjoint action of $Sp
(2n,\R)$ on $\spn$. Hence, by Proposition~\ref{prop-smoothening}, we
can assume without loss of generality that $\zeta$ is
$C^\infty$-smooth on this orbit. This yields that $F$ is
$C^\infty$-smooth on the set $\{ \xi+i\eta\in\C^{2n}\ | \ \omega
(\xi,\eta) = 1\}$. One readily checks that
\begin{equation}
\label{eqn-F-homogeneous} F(\lambda v) = |\lambda|^2 F(v)\ \forall
\lambda\in \C,
\end{equation}
which yields that $F$ is $C^\infty$-smooth on $\cW=\{ \xi +
i\eta\in \C^{2n}, \omega (\xi,\eta)> 0\}$.

\begin{lemma}
\label{lem-F-restr-is-a-Herm-qform} The restriction of $F$ on any
$L_J$, $J\in \cJ$, is a Hermitian quadratic form (for the
definitions of $L_J$ and $\cJ$ see
Section~\ref{sec-functions-whose-restr-are-qforms} above).

\end{lemma}

\begin{proof}
Given $J\in \cJ$, the space of $A\in\spn$ commuting with $J$ is a
Lie subalgebra $\gou (J)$ of $\spn$ isomorphic to the real Lie
algebra $\gou (n)$ of skew-Hermitian complex $n\times n$-matrices.
As one can easily check, $Y_{\xi, J\xi}\in \gou (J)$ for any
$\xi\in\R^{2n}$ and any $J\in \cJ$. By
Corollary~\ref{cor-Gleason}, the restriction of $\zeta$ on $\gou
(J)\cong \gou (n)$ is linear and we can write
$$F(\xi+iJ\xi) = \zeta (Y_{\xi,J\xi}) = tr (HY_{\xi,J\xi})$$
for some $H\in M_{2n} (\R)$. On the other hand,
$$tr (HY_{\xi,J\xi}) = tr (HT_{\xi,\xi}) + tr(HT_{J\xi,J\xi}) = $$
$$= tr T_{\xi, H\xi} + tr T_{J\xi, HJ\xi} = \omega (\xi, H\xi) +
\omega (J\xi, HJ\xi).$$ Thus $F(\xi+iJ\xi)=tr (HY_{\xi,J\xi})$ is a
quadratic form in $\xi$. In view of \eqref{eqn-F-homogeneous}, it
means that the restriction of $F$ on $L_J$ is a real Hermitian
quadratic form.
\end{proof}

Combining Lemma~\ref{lem-F-restr-is-a-Herm-qform} with
Proposition~\ref{prop-complex-qforms-smooth-restrictions-to-LJ}
finishes the proof of Proposition~\ref{prop-F-G-quadratic}(i).
\qed

\subsection{Evaluating a continuous Lie quasi-state on $\gogl$ on rank 1
matrices}\label{subsec-rank-1}

Here is another auxiliary proposition that we will need later.

\begin{prop}
\label{prop-qstates-gln-rank1-trace-1}

Let $\zeta$ be a continuous Lie quasi-state on the Lie algebra
$\gogl$, $n\geq 3$.  Let $\cP_1\subset\gogl$ be the set of matrices
of rank $1$. Then there exists a matrix $N\in M_n (\R)$ so that
$\zeta (A) = tr NA$ for any $A\in \cP_1$.

\end{prop}

\begin{rem}
\label{rem-qstates-gln}{\rm If $\zeta$ is a continuous Lie quasi-state
on $\gogl$ and matrices $A,B,A+B\in \gogl$ are diagonalizable over
$\R$, then
$$\zeta (A+B) = \zeta (A) + \zeta (B).$$
Indeed, any matrix diagonalizable over $\R$ is a sum of commuting
rank-1 matrices.}
\end{rem}

\bigskip
\noindent {\bf Proof of
Proposition~\ref{prop-qstates-gln-rank1-trace-1}.}

By Proposition~\ref{prop-smoothening}, we can assume without loss of
generality that $\zeta$ is smooth on any orbit of the adjoint $GL_n
(\R)$-action on $\gogl$ and, in particular, on $\cP'_1:= \{ A\in
\cP_1\ | \ tr A = 1\}$ which is such an orbit.

As before we denote by $(\cdot, \cdot)$ the Euclidean inner product
on $\R^n$. Given $\xi,\eta\in \R^n$, define an operator
$B_{\xi,\eta}\in\cP_1$ on $\R^n$ by
$$B_{\xi, \eta} (x) := (x,\eta)\xi.$$
One can easily check that if $(\xi,\eta) >0$ then $B_{\frac{
\xi}{\sqrt{(\xi,\eta)}},\frac{\eta}{\sqrt{(\xi,\eta)}}}\in\cP'_1$.
Define a function $f : \R^{2n} \to \R$ by $$f (\xi,\eta):= \zeta
(B_{\xi,\eta}).$$ Then
\begin{equation}\label{eg-f-homog}
f(\lambda \xi, \eta) = f (\xi, \lambda\eta) = \lambda f(\xi,\eta)\
\forall \lambda\in\R,\end{equation} and therefore if $(\xi,\eta)>0$,
then
$$f (\xi,\eta) = \zeta (B_{\xi,\eta}) = (\xi,\eta) \zeta
\big(B_{\frac{\xi}{\sqrt{(\xi,\eta)}},\frac{\eta}{\sqrt{(\xi,\eta)}}}\big)
= (\xi,\eta) f \big( \tfrac{\xi}{\sqrt{(\xi,\eta)}},
\tfrac{\eta}{\sqrt{(\xi,\eta)}}\big).$$ Since $\zeta$ is smooth on
$\cP'_1$, we get that $f$ is a smooth function on $$\cU:= \{
(\xi,\eta)\in \R^n\times \R^n\ |\ (\xi,\eta)>0\}\;.$$

As before, denote by $\cS^+_n (\R)$ the space of all symmetric real
positive-definite $n\times n$-matrices. For any $M\in \cS^+_n (\R)$
define an inner product $(\cdot, \cdot)_M$ on $\R^n$ by
$$(x,y)_M := (Mx,y)\;.$$
Denote by $\cS_M$ the space of all real $n\times n$-matrices
symmetric with respect to this inner product:
$$\cS_M := \{ A\in M_n (\R)\ | \ (Ax, y)_M = (x,Ay)_M\ \forall
x,y\in \R^n \}.$$ By Theorem~\ref{thm-Gleason}, for any $M\in
\cS^+_n (\R)$ there exists $T_M\in \cS_M$ such that
$$\zeta (A) = tr T_M A\ \forall A\in \cS_M.$$
Given $M\in \cS^+_n (\R)$, define
$$L_M := \{ (\xi,\eta)\in \R^n\times \R^n\ |\ \eta=M\xi\}.$$
One can easily check that for any $(\xi,\eta)\in L_M$ the operator
$B_{\xi, \eta}$ lies in $\cS_M$ and hence
$$f (\xi,\eta)= \zeta (B_{\xi,\eta}) = tr (T_M  B_{\xi,\eta}) =
tr B_{T_M \xi,\eta} = (T_M\xi, \eta).$$ Thus the restriction of $f$
on any $L_M$, $M\in \cS^+_n (\R)$, is a quadratic form. Applying
Proposition~\ref{prop-real-qforms-restrictions-LA} we get that there
exists a quadratic form  on $\R^{2n}$ which coincides with $f$ on
$\cU$. It follows from \eqref{eg-f-homog} that
\begin{equation}\label{eq-f-expl}
f(\xi,\eta) = (N\xi,\eta)\end{equation} for some matrix $N$.

Observe now that $B_{-\xi,\eta}=-B_{\xi,\eta}$ and hence
$f(-\xi,\eta) = -f(\xi,\eta)$. Furthermore, $\xi \times \eta \in
\cU$ whenever $(-\xi) \times \eta \in \cU^-$, where
$$\cU^-:= \{
(\xi,\eta)\in \R^n\times \R^n\ |\ (\xi,\eta)<0\}\;.$$ Thus
equality \eqref{eq-f-expl} holds on $\cU^-$ as well. Since $\cU
\cup \cU^-$ is dense in $\R^{2n}$, we conclude that $f(\xi,\eta) =
(N\xi,\eta)$ for all $\xi$ and $\eta$ which means that
$\zeta(B_{\xi,\eta})= tr N B_{\xi,\eta}$. This completes the
proof. \Qed

\subsection{Lie quasi-states on a Heisenberg
algebra}\label{subsec-heisenberg}

Here we will prove an auxiliary result which is, in fact, equivalent
to the claim that any (not necessarily continuous) Lie quasi-state
on a Heisenberg Lie algebra has to be linear.

\begin{prop}
\label{prop-Heisenberg}

Assume that the restriction of a function $\phi: \R^{2n}\to\R$, $n
\geq 2$ to any isotropic vector subspace of $(\R^{2n},\omega)$ is
linear (here $\omega$ is a linear symplectic form on $\R^{2n}$).
Then $\phi$ is linear.

\end{prop}

\begin{proof}

Let $e_1,\ldots, e_n, f_1,\ldots,f_n$ be the Darboux basis of
$\R^{2n}$ corresponding to the coordinate system $p_1,\ldots, p_n,
q_1,\ldots, q_n$.

Any $a\in \R^{2n}$ can be written as
$$a=a_1+\ldots+a_n, \ a_i\in {\rm Span}\, (e_i,f_i),\
i=1,\ldots,n,$$ so that
$$\omega (a_i, a_j) = 0,\ \forall i\neq j,$$
and therefore
$$\phi (a) = \phi(a_1)+\ldots+\phi(a_n).$$
Thus it suffices to prove that the restriction of $\phi$ on each
${\rm Span}\, (e_i,f_i)$, $i=1,\ldots,n$, is linear.

Restricting $\phi$ on 1-dimensional vector subspaces of $\R^{2n}$ --
which are all isotropic, of course -- we get that $\phi$ is
homogeneous.

Denote
$$k:=\phi (e_1) + \phi (f_1) - \phi (e_1 + f_1).$$
Fix $i\in\N$, $2\leq i\leq n$. We need to show that
\begin{equation}
\label{eqn-phi-a-b} \phi (a+b)=\phi (a) +\phi (b)\ \forall a,b\in
{\rm Span}\, (e_i,f_i).
\end{equation}
If $\omega (a,b) = 0$, this follows from the hypothesis of the
proposition. Otherwise, after permuting if necessary $a$ and $b$,
assume that $\omega(a,b) = -C^2 < 0$. Then the vectors $Ce_1 + a$
and $Cf_1 + b$ are $\omega$-orthogonal. Hence,
$$\phi (Ce_1) + \phi (a) + \phi (Cf_1) + \phi (b) =
\phi (Ce_1 + a) + \phi (Cf_1 + b) =
$$
$$
= \phi (Ce_1 + a + Cf_1 +b) = \phi (Ce_1 + Cf_1) + \phi (a+b).$$
Therefore, by the homogeneity of $\phi$,
$$\phi (a+b) - \phi (a) - \phi (b) = Ck, \ \forall a,b\in {\rm Span}\, (e_i,f_i).$$
Substituting $-a,-b$ instead of $a,b$ in the last equation and
using again the homogeneity of $\phi$ we get that $Ck=0$, hence
$k=0$, which yields \eqref{eqn-phi-a-b}. Thus the restriction of
$\phi$ on each ${\rm Span}\, (e_i,f_i)$, $i=2,\ldots,n$, is
linear. Switching $i\neq 1$ and $1$ we see immediately that the
restriction of $\phi$ on ${\rm Span}\, (e_1,f_1)$ is linear as
well, which finishes the proof.

\end{proof}

\subsection{Proof of
Proposition~\ref{prop-F-G-quadratic}(ii).}\label{subsec-proof-ii}

Now we will deal with the function $G$ defined in \eqref{eqn-def-G}.
Let us fix $\xi\in\R^{2n}$ and show that $G(\xi,\cdot):
\R^{2n}\to\R$ is a linear function -- since $G$ is obviously
symmetric with respect to $\xi,\eta$, this would show that $G$ is a
symmetric bilinear form on $\R^{2n}$. By
Proposition~\ref{prop-Heisenberg}, it is enough to show that
\begin{equation}
\label{eqn-G-linear-on-omega-orthogonal-vectors}
\omega (\eta_1, \eta_2) = 0\ \Longrightarrow\ G(\xi,
c_1\eta_1+c_2\eta_2) = c_1 G(\xi,\eta_1) + c_2 G(\xi, \eta_2)\
\forall c_1,c_2\in\R.
\end{equation}
Let us, indeed, assume that $\omega (\eta_1, \eta_2) = 0$. Choose
a Lagrangian subspace $L_2$ containing $\eta_1, \eta_2$.

If $\xi\in L_2$, then a direct check shows that $[Z_{\xi,\eta_1},
Z_{\xi, \eta_2}]=0$ and hence
\eqref{eqn-G-linear-on-omega-orthogonal-vectors} follows by the
definition of a Lie quasi-state.

If $\xi\notin L_2$, then, as one can easily check, there exists a
Lagrangian subspace $L_1$ transversal to $L_2$ and containing
$\xi$. Define a Lie subalgebra
$$\cR (L_1, L_2):= \{ A\in \spn\ | \ AL_1\subset L_1,\ AL_2\subset
L_2\}$$ of $\spn$. An easy check using the linear Darboux theorem
shows that the mapping $A \mapsto A|_{L_1}$ establishes a Lie
algebras isomorphism between $\cR (L_1, L_2)$ and
${\mathfrak{gl}\,(L_1)} \approx \gogl$. Furthermore, any
transformation $Z_{\xi,v}$ with $v \in L_1$ lies in $\cR (L_1,
L_2)$. Its image under the above isomorphism has rank $1$ provided
$v \neq 0$. Applying
Proposition~\ref{prop-qstates-gln-rank1-trace-1} to the elements
$Z_{\xi,\eta_1}, Z_{\xi,\eta_2}$ and $Z_{\xi, \eta_1 + \eta_2} =
Z_{\xi,\eta_1} + Z_{\xi,\eta_2}$ of $\cR (L_1, L_2)$ yields
\eqref{eqn-G-linear-on-omega-orthogonal-vectors}.

This finishes the proof of Proposition~\ref{prop-F-G-quadratic} and
hence of Theorem~\ref{thm-main}.\Qed

\section{Discontinuous Lie quasi-states}
\label{sec-discont-Lie-qs}

In this section we prove Theorem~\ref{thm-discont}. We start from
the following general observation.

\begin{prop}
\label{prop-discont-Lie-qs-general} Assume $L\subset \gog$ is an
abelian subalgebra of a Lie algebra $\gog$ and $L_0\subset L$ is a
vector subspace so that
$$[x,v]=0, x\in \gog, v\in L\setminus
L_0\ \Longrightarrow\ x\in L.$$ Let $\alpha: L\to\R$ be a linear
functional on $L$ such that $\alpha\not\equiv 0$ and
$\left.\alpha\right|_{L_0}\equiv 0$. Define a functional $\zeta:
\gog\to\R$ as
$$\zeta (x) =
\begin{cases}
  0,\ \hbox{\rm if}\ \, x\notin L, \\
  \alpha (x),\ \hbox{\rm if}\ \, x\in L. \\
\end{cases}$$
Then $\zeta$ is a Lie quasi-state.

\end{prop}

\medskip
\noindent \begin{proof} Assume that $[x,y]=0$. We have to show that
\begin{equation}\label{vsp-prop-disc-10} \zeta(x+y)=\zeta(x)+\zeta(y)\;.\end{equation}
Vectors $x,y,x+y$ pairwise commute. If at least one of them does not lie in $L$, two others
must lie in $(\gog \setminus L)\cup L_0$. Thus $\zeta$ vanishes on each of these vectors and
so \eqref{vsp-prop-disc-10} holds. If $x,y \in L$ equation \eqref{vsp-prop-disc-10} follows from the linearity of $\alpha$.\end{proof}

\medskip

In the case of $\gog=\spn$ one can construct $L,L_0$ as in
Proposition~\ref{prop-discont-Lie-qs-general} in the following
way.

\begin{lemma}
\label{lem-odd-power-skew-sympl-matrix} Let $A\in\spn$ and let
$p(t)$ be a real polynomial. Then $p(A)\in\spn$ if and only if
$p(t)$ includes only the odd powers of $t$ (i.e. $p$ is an odd
function of $t$).
\end{lemma}

\medskip
\noindent This is an immediate consequence of the fact that
$A\in\spn$ if and only if $AJ=-JA^T$, where $J$ is the
standard complex structure on $\R^{2n}$.

Take any matrix $A\in \spn$  whose Jordan form is a $2n\times 2n$
Jordan block with the eigenvalue $0$ (the existence of such an $A$
is well-known -- see e.g. \cite{Williamson}, cf.
\cite{Arnold-math-methods}, \cite{Hormander}). Define $L_0, L$ as
follows:
$$L:=\{ a_1 A + a_3 A^3 +\ldots + a_{2n-1} A^{2n-1},\ a_1, a_3,\ldots,
a_{2n-1}\in\R \},$$
$$L_0:= \{ a_3 A^3 +\ldots + a_{2n-1} A^{2n-1},a_3,\ldots,
a_{2n-1}\in\R \}\subset L.$$ Let us show that subspaces $L_0$ and
$L$ satisfy the hypothesis of
Proposition~\ref{prop-discont-Lie-qs-general}.

\begin{lemma}
\label{lem-commutes-with-powers-of-A}
$$[x,v]= 0, v\in L\setminus L_0\ \Longrightarrow [x,A^{2m-1}] = 0\
\forall m=1,\ldots, n.$$
\end{lemma}

\begin{proof} Without loss of generality,
$$ v = A+a_3 A^3 + ... + a_{2n-1}A^{2n-1}\;.$$
We prove the statement of the lemma using an inverse induction by
$m$ (from $m=n$ to $m=1$).

For $m=n$ note that $v^{2n-1} = A^{2n-1}$ (we use here that
$A^{2n}=0)$, and hence $[x,v^{2n-1}]=[x,A^{2n-1}]=0$.

Assume now that $[x,A^{2j-1}]=0$ for all $j = m+1,...,n$. Put
$$B=A+\sum_{i=2}^{m} a_{2i-1}A^{2i-1}\;.$$ The inductive
assumption together with $[x,v]=0$ yields $[x,B]=0$. But
$$B^{2m-1}= A^{2m-1} + \sum_{i \geq 2m+1}  c_i A^i$$ for some
coefficients $c_i$. Together with the inductive assumption this
yields $[x,A^{2m-1}]= [x,B^{2m-1}]=0$, as required.
\end{proof}

Finally, we recall without proof a standard fact from the linear
algebra.

\begin{lemma}
\label{lem-commutes-with-Jordan-block} Assume that a real square
matrix $\Delta$ is a Jordan block with the eigenvalue zero. Then any
matrix commuting with $\Delta$ is a polynomial of $\Delta$.
\end{lemma}

\medskip
\noindent {\bf Proof of Theorem~\ref{thm-discont}:} Assume that
$v\in L\setminus L_0$ and $[x,v]=0$. Then, by
Lemma~\ref{lem-commutes-with-powers-of-A}, $[x,A] = 0$. Hence, by
Lemma~\ref{lem-commutes-with-Jordan-block}, $x$ is a polynomial of
$A$, which, by Lemma~\ref{lem-odd-power-skew-sympl-matrix},
belongs to $L$. This shows that subspaces $L_0$ and $L$ satisfy
the hypothesis of Proposition~\ref{prop-discont-Lie-qs-general}.
Applying this proposition and varying the matrix $A$ together with a functional $\alpha$,
we get a continuum of linearly independent
discontinuous Lie quasi-states on $\spn$. All of them are bounded
in a neighborhood of zero. \Qed

\section{$Ad$-invariant Lie quasi-states}
\label{sec-Ad-inv-Lie-qs}

In this section we will discuss Lie quasi-states invariant under
the adjoint action of a Lie group.

Note that the adjoint actions on $\gog$ of different connected Lie
groups with the same Lie algebra $\gog$ have the same orbits because
such groups are locally isomorphic (see e.g. \cite{Helgason}, p.
109) so that the local isomorphism intertwines their adjoint
actions, and any connected Lie group is generated by a neighborhood
of the identity. We say that a function on the Lie algebra $\gog$ is
{\it $Ad$-invariant} if it is constant on any orbit of the adjoint
action on $\gog$ of any Lie group with the Lie algebra $\gog$.

Note also that if $G_1$ is a (closed) Lie subgroup of a Lie group
$G_2$ and $\gog_1\subset\gog_2$ are the corresponding Lie
algebras, then any $Ad$-invariant function on $\gog_2$ (that is,
invariant with respect to the adjoint action of $G_2$) restricts
to an $Ad$-invariant function on $\gog_1$ (that is, invariant with
respect to the adjoint action of $G_1$).

Recall that a real Lie algebra is called {\it compact} if it is
the Lie algebra of some compact real Lie group. A Lie group $G$
will be called {\it Hermitian} (see e.g. \cite{K}) if
\begin{itemize}
\item{} $G$ is connected and non-compact;

\item{} the Lie algebra of $G$ is simple;

\item{} the associated homogeneous space $G/K$, where $K$ is the
maximal compact subgroup of $G$, has a complex-manifold structure
and $G$ acts on it by holomorphic transformations.

\end{itemize}
There is a complete classification of the Lie algebras of
Hermitian Lie groups (see e.g. \cite{K}). In particular, $Sp\,
(2n,\R)$ is a Hermitian Lie group.

We will now classify $Ad$-invariant Lie quasi-states on compact and
Lie algebras and the Lie algebras of Hermitian Lie groups. Let us
emphasize that a priori these quasi-states are not assumed to be
continuous.

\begin{thm}
\label{thm-no-Lie-qs-on-compact-Lie algebras} Any  $Ad$-invariant
Lie quasi-state on any compact Lie algebra is linear. If the
compact Lie algebra has a trivial center (in particular, if it is
simple), any $Ad$-invariant Lie quasi-state on it vanishes
identically.

\end{thm}

\begin{thm}
\label{thm-Lie-qs-Hermitian-Lie-algebras} Let $G$ be a simply
connected Hermitian Lie group and let $\gog$ be its Lie algebra. Let
$\mu$ be the unique (up to a multiplicative constant) non-trivial
homogeneous quasi-morphism on $G$ and let $\xi: \gog\to\R$ be its
directional derivative:
$$\xi(x):= \mu(\exp(x)) \;\;\forall x \in \gt{g}\;.$$
Then any  $Ad$-invariant Lie quasi-state $\zeta$ on $\gog$ is
proportional to $\xi$.
\end{thm}

\medskip
\noindent We refer to \cite{Shtern} (cf. \cite{BS-H}) for the
uniqueness (up to a multiplicative constant) of a non-trivial
homogeneous quasi-morphism on $G$.

\bigskip
\noindent {\bf Proof of Theorem~\ref{thm-no-Lie-qs-on-compact-Lie
algebras}:} Any compact Lie algebra can represented as a direct sum
of an abelian Lie algebra and a number of compact simple Lie
algebras (see e.g. \cite{Helgason}, p. 132). Any Lie quasi-state on
an abelian Lie algebra is linear. This shows that the first claim of
the theorem follows from the second one.

Let us show that any $Ad$-invariant Lie quasi-state $\zeta$ on a
compact simple Lie algebra $\gog$ vanishes identically -- this
would immediately imply the second claim of the theorem. Denote by
$G$ a compact connected Lie group whose Lie algebra is $\gog$. Any
element of $\gog$ lies in a Cartan subalgebra of $\gog$, that is
the abelian subalgebra which is the Lie algebra of a maximal torus
in $G$, and any two Cartan subalgebras are mapped into each other
by the adjoint action of $G$ (see e.g. \cite{DK}, p.152). Thus it
suffices to show that $\zeta$ vanishes on a Cartan subalgebra
$\goh\subset \gog$. Since $\goh$ is abelian, the restriction of
$\zeta$ to $\goh$ is linear. Since $\zeta$ is $Ad$-invariant, this
linear function on $\goh$ is invariant under the actions of
$Ad_g$, $g\in G$, that preserve $\goh$, that is under the action
of the Weyl group $W$ of $G$ on $\goh$. Since the action of $W$ on
the simple Lie algebra $\goh$ has only trivial invariant subspaces
(see e.g. \cite{DK}, p. 172 and p.251), we have
$\left.\zeta\right|_\gog \equiv 0$, as required. \Qed

\bigskip
\noindent {\bf Proof of
Theorem~\ref{thm-Lie-qs-Hermitian-Lie-algebras}:}

\medskip
\noindent 1) The structure of the Hermitian  Lie group on $G$
gives rise to the following features of the Cartan decomposition
$\gt{g} = \gt{t} + \gt{p}$, where $\gt{t}$ is the Lie algebra of a
maximal compact subgroup $K$ of $G$: The center $\gt{c}$ of
$\gt{t}$ is one-dimensional and contains a preferred element, say
$J$, so that $\ad_J$ preserves $\gt{p}$ and acts on $\gt{p}$ as a
complex structure (see e.g. \cite{K}, Theorem 7.117  and p. 513).
We shall normalize $\xi$ and $\zeta$ by $\xi(J)=\zeta(J)=1$.

\medskip
\noindent 2) Let us check that any $Ad$-invariant Lie quasi-state
must vanish on $\gt{p}$ by using a trick by Ben Simon and Hartnick
\cite{BS-H}: they noticed that since $(\ad_J)^2 = -\bf{1}$, one
has $\exp (\pi \cdot \ad_J) = -\bf{1}$ on $\gt{p}$ (here $\pi$ is
the number $\pi=3.14\ldots$). It follows that $\zeta(x) =
-\zeta(x)$ for all $x \in \gt{p}$ which yields the claim.

\medskip
\noindent 3) Every element $x \in \gt{g}$ can be written as $s+n$,
where $s$ is semi-simple, $n$ is nilpotent and $[s,n]=0$. By
Jacobson-Morozov theorem (\cite[p.620]{K}, pass to the
complexification and use that $n$ is a real nilpotent) , there
exists $y \in \gt{g}$ with $[y,n] = n$. Therefore, setting
$f=\exp(y)$ we get (passing to the series for the exponent) that
$Ad_f(n) = e\cdot n$, where $e$ is the number $e=2.71\ldots$. This
yields
$$\zeta(n) =\zeta(e\cdot n) = e \cdot \zeta(n)\;.$$ Therefore $\zeta
(n)=0$ and so $\zeta(x) = \zeta(s)$. Thus it suffices to check
that $\zeta=\xi$ on the semi-simple elements.

\medskip
\noindent 4) Every semi-simple element of $\gt{g}$ lies in some
Cartan subalgebra. Every Cartan subalgebra of $\gt{g}$ is
conjugate to a Cartan subalgebra of the form   $\gt{a} + \gt{h}$
where $\gt{a}\subset \gt{p}$ and $\gt{h} \subset \gt{t}$
(\cite[Proposition 6.59]{K}). Thus every semi-simple element of
$\gt{g}$ is conjugate to $a+h$ with $a \in \gt{a}, h\in\gt{h}$ and
$[a,h]=0$. Since $\zeta$ and $\xi$ vanish on $\gt{a}$ (step 2), it
suffices to show that $\zeta=\xi$ on $\gt{h}$. In fact, we shall
show that any $Ad$-invariant Lie quasi-state $\zeta: \gog\to\R$
vanishes on the algebra $\gt{t}':= [\gt{t},\gt{t}]$, which would
yield the desired result in view of the decomposition $\gt{t}=
\gt{t}' \oplus \gt{c}$, where $\gt{c}$ is the center of $\gt{t}$.
Indeed, $\gt{t}'$ is a compact Lie algebra with a trivial center
(see e.g. \cite[p.513]{K}) and the restriction of
$\left.\zeta\right|_{\gt{t}'}$ is an $Ad$-invariant Lie
quasi-state on $\gt{t}'$. Therefore, by
Theorem~\ref{thm-no-Lie-qs-on-compact-Lie algebras}, $\zeta$
vanishes identically on $\gt{t}'$, which finishes the proof.\Qed

\bigskip \noindent {\bf Acknowledgements.} We thank Semyon Alesker
for explaining to us the smoothening procedure used in
Section~\ref{subsec-smooth} and Frol Zapolsky for correcting an
inaccuracy in the first version of the paper. The early results of
the paper were presented at seminars in Tel Aviv and Neuchatel. We
thank the organizers, Vitali Milman and Felix Schlenk, for this
opportunity.

\bibliographystyle{alpha}

\bigskip

\noindent
\begin{tabular}{ll}
Michael Entov & Leonid Polterovich\\
Department of Mathematics & School of Mathematical Sciences\\
Technion - Israel Inst. of Technology & Tel Aviv University\\
Haifa 32000, Israel & Tel Aviv 69978, Israel\\
entov@math.technion.ac.il & polterov@post.tau.ac.il\\
\end{tabular}

\end{document}